\RequirePackage{fix-cm}
\documentclass[smallextended]{svjour3}       
\smartqed  
\usepackage[margin=1.2in]{geometry}
\usepackage{graphicx}
\usepackage{xcolor}  

\usepackage[utf8]{inputenc}  
\usepackage{enumerate}

\usepackage{amsmath}  
\usepackage{amsfonts}  
\usepackage{booktabs} 
\usepackage{multirow} 

\usepackage{doi}  

\newcommand\Real{\mathbb{R}}  
\newcommand\Lag{\mathcal{L}}  
\DeclareMathOperator*{\minimize}{minimize}

\begin{document}

\title{OpenSQP: A Reconfigurable Open-Source SQP Algorithm in Python for Nonlinear Optimization
}


\author{Anugrah Jo Joshy$^1$        \and
        John T. Hwang$^1$  
}


\institute{Anugrah Jo Joshy \at
              \email{ajoshy@ucsd.edu}           
}

\date{Received: date / Accepted: date}

\maketitle

\footnotetext[1]{Department of Mechanical and Aerospace Engineering, 
    University of California San Diego, La Jolla, CA 92093, USA}

\begin{abstract}
Sequential quadratic programming (SQP) methods have been remarkably 
successful in solving a broad range of nonlinear optimization problems.
These methods iteratively construct
and solve quadratic programming (QP) subproblems
to compute directions that converge to a local minimum.
While numerous open-source and commercial SQP algorithms are available, 
their implementations lack the transparency  and modularity necessary 
to adapt and fine-tune them for specific applications
or to swap out different modules to create a new optimizer.
To address this gap, we present OpenSQP, a modular and reconfigurable 
SQP algorithm implemented in Python
that achieves robust performance comparable to leading algorithms.
We implement OpenSQP in a manner that allows users to easily modify 
or replace components such as merit functions, 
line search procedures, Hessian approximations, and 
QP solvers.
This flexibility enables the creation of
tailored variants of the algorithm for specific needs.
To demonstrate reliability, we present numerical results using the standard 
configuration of OpenSQP that employs a smooth augmented 
Lagrangian merit function for the line search and 
a quasi-Newton BFGS method for approximating the Hessians.
We benchmark this configuration
on a comprehensive set of problems from 
the CUTEst test suite.
The results demonstrate performance that is competitive with
proven nonlinear optimization algorithms such as 
SLSQP, SNOPT, and IPOPT.



\end{abstract}

\section{Introduction}
\label{sec:intro}

Nonlinear programming (NLP) refers to the 
solution of
optimization problems characterized by 
nonlinear, continuously differentiable objective and 
constraint functions.
Such problems are frequently encountered in various contexts,
including engineering design, financial analysis, and
economic modeling.
In engineering, one of the most common applications of NLP 
is in the optimal design and control of complex systems
such as aircraft, robotic systems, and chemical processing plants.

Over the past few decades, sequential quadratic programming (SQP) and 
interior point (IP) methods have established themselves as the most efficient 
approaches for solving general nonlinear optimization problems.
State-of-the-art SQP methods, such as SNOPT \cite{gill2005snopt},
have been shown to be effective for problems with up to
tens of thousands of variables and constraints \cite{hwang2014large}.
However, their efficiency deteriorates as the problem size increases, 
primarily due to the challenges associated with identifying 
the set of active constraints at the solution.
IP methods, on the other hand, overcome this limitation 
by solving a sequence of equality-constrained barrier subproblems, 
where all inequality constraints are transformed into equalities 
through the introduction of slack variables.
IP methods have proven effective in solving large-scale problems with 
millions of variables \cite{wachter2002interior}.

Although IP methods are generally considered more
reliable for solving nonlinear problems when no prior information
about the solution's properties is available, 
numerical studies suggest otherwise \cite{gill2015performance}.
A rigorous comparison of IP and SQP methods remains 
challenging, as most software packages based on IP methods 
utilize second derivatives, whereas those based on SQP typically do not.
Comparative studies relying solely on first derivatives indicate that,
on average, SQP methods demonstrate greater efficiency and 
higher success rates compared to IP methods \cite{gill2015performance}.
Nevertheless, IP methods may still be better suited for certain problems,
even when limited to first-derivative information.
For example, quasi-Newton interior point methods frequently 
outperform their SQP counterparts when solving trajectory 
optimization problems starting from poor initial guesses \cite{kaneko2025simultaneous}.

In practice, obtaining even first derivatives presents 
significant challenges in many problems. 
For instance, in multidisciplinary design optimization (MDO), 
the computation of first derivatives for complex, large-scale models 
has only recently become feasible due to advancements 
in modeling frameworks that integrate automatic differentiation methods
\cite{gray2019openmdao,gandarillas2024graph}. 
Similarly, there are cases where computation of second derivatives is 
costly or impractical.
In such scenarios, SQP methods naturally 
become the default choice.

In addition to their superior performance on problems with 
only first-order information, SQP methods are also
better at leveraging good initial guesses, certifying infeasible constraints,
exploiting linear constraints, and achieving high accuracy solutions.
The first three properties are particularly crucial in applications that 
involve solving a sequence of related NLP problems, 
e.g., mixed-integer nonlinear programming (MINLP),
successive optimization with refinement of an underlying discretization,
and system-level design optimization during the conceptual phase
where subsystem models undergo continuous revisions.
For these reasons, SQP remains an active area of research.

One of the drawbacks of conventional SQP methods, compared to IP methods,
is their inability to readily utilize Hessians when they are available.
SQP methods typically rely on positive definite quasi-Newton Hessian approximations
to circumvent the difficulties arising from nonconvex QP subproblems generated 
by the use of exact Hessians.
Recent developments in SQP theory focus on effectively incorporating 
second derivative information into QP subproblems
through the convexification of Hessians \cite{gill2014convexification}.

Another drawback of SQP is its limited compatibility with general-purpose
linear solvers optimized for modern computer architectures.
QP subproblems are often solved using active-set
methods for efficiency, but this results in linear systems 
that change structure between QP iterations,
hindering the effective use of state-of-the-art linear algebra software.
In comparison, the Newton equations in IP methods always include 
the gradients of all inequality constraints, preserving the dimensions and
sparsity structure of the linear systems throughout the optimization.
Recent advancements in SQP software emphasize leveraging 
off-the-shelf linear solvers.

The advent of modern programming languages such as
Python and Julia has transformed scientific computing.
However, performance-focused optimization algorithms are still largely 
written in low-level languages such as C or Fortran
to retain the computational efficiency of compiled languages.
Many widely used algorithms were originally developed 
decades ago, and the substantial effort required to modernize them 
has led to their continued development in the same low-level languages.
While some algorithms have been implemented in high-level languages, 
they are primarily educational in nature and often lack the 
efficiency, robustness, or versatility needed for general-purpose optimization.
The availability of transparent and reliable algorithms is crucial
for research in nonlinear programming, as practical implementation
and testing often serve as validation for novel algorithms, especially since 
the theoretical characterization of algorithm performance is challenging,
even in idealized settings.
Aside from the \textit{trust-constr} algorithm in the 
SciPy \cite{virtanen2020scipy} library, 
few transparent and robust NLP algorithms are available.
This lack of accessible implementations poses a barrier for 
engineers and researchers seeking to modify these algorithms 
for their specific use cases.

To address this challenge, we present OpenSQP, a modular and reconfigurable SQP algorithm
implemented in Python that delivers robust performance comparable to 
state-of-the-art algorithms.
Our implementation allows users to easily modify 
or replace key components such as merit functions, line search procedures,
Hessian approximations, and QP solvers. 
This flexibility facilitates the development of custom SQP optimizers 
that are fine-tuned to meet the needs of specific applications.
OpenSQP also supports researchers developing new SQP algorithms
by providing an accessible starting point, 
even if they are novice programmers.
To demonstrate the reliability of our implementation, we present numerical results 
using the standard configuration of OpenSQP that employs a smooth 
augmented Lagrangian merit function for the line search and
a quasi-Newton BFGS method for approximating the Hessians.

For broader adoption and ease of distribution, a new optimization
algorithm is best implemented within an existing optimizer library.
These libraries also provide utilities such as tools for
recording, warm-starting, visualizing, and post-processing, 
which new algorithms can inherit without additional implementation effort.
Several Python libraries exist for nonlinear programming, 
including pyOpt \cite{perez2012pyopt}, pyOptSparse \cite{wu2020pyoptsparse}, 
SciPy’s \texttt{optimize} module \cite{virtanen2020scipy}, 
and modOpt \cite{joshy2024modopt}.
Among these, only modOpt is explicitly designed to support the 
development of modular, transparent, and reconfigurable optimizers.
Unlike the others, which primarily cater to advanced developers, 
modOpt provides dedicated tools for optimizer development, 
including an interface to the CUTEst \cite{gould2015cutest,fowkes2022pycutest} 
test problem collection and utilities for performance profiling.
For these reasons, we implement and distribute OpenSQP
as part of the modOpt library.
As a result, the algorithm is available as open-source on modOpt's 
GitHub repository, 
with its documentation provided through modOpt's documentation.

The remainder of the paper is structured as follows.
Section \ref{sec:bg} introduces fundamental concepts of nonlinear programming
and sequential quadratic programming, along with the most commonly employed 
approaches for each.
Section \ref{sec:meth} details the algorithmic components of
the standard configuration of OpenSQP.
Section \ref{sec:nr} evaluates the reliability and efficiency of 
the standard configuration on various benchmark problems.
Section \ref{sec:conclusion} concludes the paper and 
outlines potential directions for future work.

\section{Background}
\label{sec:bg}

\subsection{The Problem}
\label{sec:prob}
The algorithm we discuss in this paper applies to the nonlinear
programming problem
\begin{equation}
    \minimize_{x\in\Real^n} f(x) \quad \text{subject to} \quad c(x)\geq0,
    \label{eq:np}
    \tag{NP}
\end{equation}
where the objective function $f:\Real^n\to\Real$ and 
constraint functions in $c:\Real^n\to\Real^m$
are assumed to be continuously differentiable with respect
to the decision variables $x$.
While these functions may be nonlinear in $x$, we assume
their first derivatives are available, 
but make no further assumptions about $f$ and $c$.

For notational convenience, we transform all constraints into
inequalities with zero lower bounds.
Although this constraint format makes no distinction between 
variable bounds, linear constraints, and equality constraints,
the methods we describe can be easily adapted to take advantage
of these special types of constraints.

\subsection{Notation}

We denote the gradient of the objective function by $g(x)$
and the Jacobian of the constraint functions by $J(x)$.
The Lagrangian function for the NLP problem (\ref{eq:np})
is defined as $\Lag(x, \lambda)= f(x) - \lambda^T c(x)$, where 
$\lambda \in \Real^{m}$ corresponds to the vector of Lagrange multipliers
associated with the inequality constraints $c(x) \geq 0$.
The Hessian of the Lagrangian with respect to $x$ is expressed as
$H(x, \lambda) = \nabla_{xx}^2 \Lag(x, \lambda) = 
\nabla^2 f(x) - \sum_{i=1}^{m} \lambda_i \nabla^2  c_i(x)$.
We use $s$ to denote the vector of slack variables, and ($x^*$, $\lambda^*$, $s^*$)
represents the primal, dual, and slack variables 
at a local solution of (\ref{eq:np}).
The solution estimates produced by an optimization algorithm
after the $k$th iteration are denoted by ($x_k$, $\lambda_k$, $s_k$).
Additionally, we use $f_k$, $g_k$, $c_k$, and $J_k$ to denote the
optimization functions and their derivatives evaluated at $x_k$.

\subsection{Optimality conditions}
\label{sec:optimality}

Optimization algorithms generate a sequence of decision variable 
iterates $\{x_k\}$ that converges to a local solution $x^*$ of the problem.
A point $x^*$ is a local solution if it satisfies certain optimality conditions.
The optimality conditions can certify a point as a solution 
only if certain assumptions about the constraints, known as constraint qualifications,
hold at that point.
In all discussions in this paper, we assume
the linear independence constraint qualification (LICQ),
which holds at a point $x$ if the gradients of the active 
inequality constraints are linearly independent.
Active constraints are those that are satisfied as equalities at a given point.

For the inequality-constrained problem (\ref{eq:np}),
the first-order optimality is characterized by
\begin{equation}
    c(x^*)\geq0, \quad \lambda^* \geq 0, \quad  \nabla_x \Lag(x^*, \lambda^*)  = 0, \quad c(x^*)^T \lambda^* = 0.
    \label{eq:kkt}
    \tag{KKT}
\end{equation}
These first-order necessary optimality conditions, also 
known as the \textit{Karush--Kuhn--Tucker (KKT)} conditions, are satisfied
at any local solution where LICQ holds.
Points satisfying the KKT conditions are referred to as 
KKT points, critical points, or stationary points.
While first-order optimality conditions are necessary, they are not sufficient to
qualify a point as a solution.

To validate a point as a solution, one must also ensure the satisfaction of
second-order sufficient conditions.
Several second-order sufficient conditions have been formulated
under assumptions of varying levels of strictness \cite{gill2020computational}.
A strong version of the sufficient optimality conditions for an
isolated local solution is given below.
A point $x^*$ is an isolated local solution if there exists a
sufficiently small neighborhood around $x^*$ where it is the only solution.

\begin{theorem}
A point $x^*$ is an isolated local solution of (\ref{eq:np}) if
\begin{enumerate}[i)]
    \item the LICQ holds at $x^*$,
    \item $x^*$ is a KKT point,
    \item strict complementarity holds at $x^*$, i.e.,
    the multiplier $\lambda^*_i > 0$ for all $i$ with $c_i(x^*)=0$, and
    \item $p^T H(x^*, \lambda^*) p > 0$ for all $p \neq 0$ such that
    $J_a(x^*)p = 0$, where $J_a(x^*)$ represents the Jacobian of the active constraints.
\end{enumerate}

\end{theorem}

\subsection{NLP methods}
\label{sec:nlp_meth}

Many general NLP methods have been proposed during the latter half of the past century.
These methods can be broadly classified into sequential unconstrained methods,
sequential quadratic programming (SQP) methods, and interior point (IP) methods,
each differing in how they handle the nonlinear inequality constraints.
Early approaches focused on minimizing a sequence of unconstrained problems,
as techniques for solving unconstrained problems were the most developed at the time.
SQP methods gained prominence in the late 1970s, 
followed by IP methods a decade later.
Today, SQP and IP methods are considered state-of-the-art in nonlinear programming.

Sequential unconstrained methods can be further categorized as
penalty methods or augmented Lagrangian methods, depending on how 
their objectives are formulated.
They solve a sequence of unconstrained subproblems 
with modified objectives that include additional terms
for penalizing constraint violations.
In contrast, SQP methods solve a sequence of 
constrained quadratic programming (QP) subproblems,
each approximating the original problem at the current solution estimate.
IP methods, on the other hand, solve a sequence of
equality-constrained barrier subproblems based on certain perturbed
optimality conditions of the original problem.
For more details on these classes of methods, 
see \cite{nocedal1999numerical} or \cite{martins2021engineering}.

Several successful software packages were developed based on
a wide variety of strategies.
Among the most popular augmented Lagrangian methods are 
MINOS \cite{murtagh1978large,murtagh1982projected,murtagh1998minos}
and LANCELOT \cite{conn1992lancelot,conn1996numerical,gould2003galahad}.
MINOS follows a linearly constrained Lagrangian (LCL) approach, which typically
requires fewer subproblem solutions but more function and derivative evaluations
than SQP or IP methods.
As a reduced-Hessian method, MINOS is well-suited for large-scale problems 
with moderate degrees of freedom (e.g., up to 1000).
LANCELOT, on the other hand, follows a bound-constrained Lagrangian (BCL) approach,
and can efficiently handle problems with tens of thousands of degrees of freedom 
\cite{bongartz1997numerical1,bongartz1997numerical2}.
Additionally, it can utilize second derivatives when available.

CONOPT \cite{drud1985conopt,drud1994conopt} is a generalized reduced gradient (GRG) 
method that has proven effective for problems with up to 500 degrees of freedom.
However, it generally demands more function and derivative evaluations 
than modern SQP or IP methods.
Although CONOPT is a reduced-Hessian method, it does not 
fall under any of the three categories described earlier.
Recent versions of CONOPT can take advantage of second derivatives.

NLP software based on SQP and IP methods have been more successful and
widely adopted due to their superior convergence properties and scalability.
NPSOL \cite{gill1986some,gill2001npsol}, NLPQL \cite{schittkowski1986nlpql},
and SLSQP \cite{kraft1988software},
were among the earliest widely used SQP packages.
The first two were proprietary software using an active-set approach,
capable of solving problems with several thousand variables and constraints
while requiring fewer function evaluations than competing methods at the time. 
NPSOL was later succeeded by SNOPT \cite{gill2005snopt}, and NLPQL by 
NLPQLP \cite{schittkowski2008sequential,schittkowski2011robust,schittkowski2015nlpqlp},
both of which can efficiently handle problems with tens of thousands of variables 
and constraints.
SLSQP, although open-source, 
has seen limited updates since its initial development in the 1980s.
However, it remains widely used due to its open-source availability and
easy accessibility through the SciPy library.
 
The interior point revolution began in the mid-1980s, with 
mature software packages competitive with state-of-the-art SQP methods 
emerging toward the end of the twentieth century.
Modern interior point solvers,
such as IPOPT \cite{wachter2006implementation}, LOQO \cite{vanderbei1999loqo}, and
KNITRO \cite{byrd1999knitro,byrd2000knitro}, have been available since the early 2000s.
These packages demonstrated the ability to solve larger-scale problems 
than those tractable by the most advanced SQP methods.

A common theme among the leading algorithms discussed above is that
their development began decades ago, and
they are implemented in low-level languages, primarily C or Fortran.
This underscores the inherent challenges in developing competitive NLP algorithms.
The only exception to this trend is SciPy's \textit{trust-constr} algorithm,
a trust-region interior point method
\cite{lalee1998implementation,byrd1999knitro,conn2000trust},
introduced in 2018.
Unlike other solvers, \textit{trust-constr} is fully implemented in Python, 
offering transparency and ease of modification.
While not as powerful as its counterparts, its transparent and open-source nature 
facilitates contributions and research in IP algorithms.
The reconfigurable SQP algorithm presented in this paper
offers a viable alternative to \textit{trust-constr}.
However, unlike \textit{trust-constr}, 
our algorithm demonstrates competitiveness with state-of-the-art solvers,
as shown in Section \ref{sec:nr}.

\subsection{Sequential Quadratic Programming}
\label{sec:sqp}

Sequential quadratic programming was first introduced by 
Wilson \cite{wilson1963simplicial} in the context of convex programming.
Subsequent work by Biggs \cite{biggs1972constrained}, 
Han \cite{han1976superlinearly,han1977globally}, 
and Powell \cite{powell1977fast,powell1978algorithms,powell1978convergence} 
extended this idea to practical SQP algorithms for general nonlinear optimization.
The following provides a concise overview of sequential quadratic programming.

SQP methods solve the general inequality-constrained problem (\ref{eq:np}) by
replacing it with a sequence of quadratic programming (QP) subproblems,
each approximating the original problem at the current iterate $(x_k,\lambda_k)$
to determine the direction leading to the next iterate $(x_{k+1},\lambda_{k+1})$.
The outer iterations that generate the sequence $(x_k,\lambda_k)$
are typically referred to as \textit{major} iterations,
whereas the iterations performed within the QP solver
are referred to as \textit{minor} iterations.
The QP subproblem at the $(k+1)$th iteration is given by
\begin{equation}
\begin{split}
    \minimize_{p_k\in\Real^n} & \quad f(x_k) + g(x_k)^T p_k + 
            \frac{1}{2} p_k^T H(x_k,\lambda_k) p_k \\
    \text{subject to} & \quad c(x_k) + J(x_k)p_k \geq 0,
\end{split}
\label{eq:qp}
\tag{QP}
\end{equation}
where the objective is a quadratic approximation of the Lagrangian and 
the constraints are linearizations of the original nonlinear constraints.
A line-search SQP method minimizes a merit function along the direction
provided by the solution $(p_k, \widehat \lambda_k)$ of (\ref{eq:qp})
to compute the step size $\alpha_k$.
The new iterates are then computed as: 
\begin{subequations}
\begin{align}
    x_{k+1}  & = x_k + \alpha_k p_k,  \\
    \lambda_{k+1} & = \lambda_k + \alpha_k (\widehat \lambda_k - \lambda_k).
    \label{eq:sqp-ls}
\end{align}
\end{subequations}
While we primarily focus on line search strategies to ensure global convergence, 
we note that trust-region methods are also widely used in practice.

\subsection{SQP approaches}
\label{sec:sqp_app}

Although all SQP methods share the fundamental approach described above,
there exists a wide variety of SQP algorithms based on choices of
QP solver, globalization strategy, and Hessian approximation.
The most efficient SQP implementations employ an active-set strategy, 
wherein each QP subproblem is solved through 
a sequence of equality-constrained QP problems that iteratively estimate 
the active set of constraints at the QP solution.
Each QP subproblem is initialized with the active set identified at the solution
of the previous QP subproblem.
As the major iterations progress, the active set estimates stabilize and 
the QP subproblems (QPs) require fewer iterations to solve,
with the final QPs often solved in a single minor iteration. 
This warm-start capability makes active-set methods particularly attractive for 
SQP implementations, where good initial estimates of the active set become 
available after the first few major iterations.

While active-set methods exhibit polynomial-time complexity for most practical problems,
their worst-case complexity is exponential in the problem size due to 
the combinatorial nature of identifying correct active sets.
As a result, early QP subproblems can become challenging to solve for large-scale 
problems, especially when the initial iterates are far from the optimum 
and the active set changes substantially between major iterations.
Moreover, computing exact solutions with the correct active set 
for early QP subproblems is wasteful when the active-set 
is still evolving significantly.
To mitigate this, several strategies based on inexact QP solutions have been proposed.
For instance, it has been shown that incomplete QP solutions that terminate at
a stationary point can still yield a valid search direction 
and ensure global convergence \cite{murray1995sequential}.

Primal active-set QP solvers typically operate in two phases: 
Phase 1 computes a feasible point, while Phase 2 minimizes 
the objective within the feasible region.
For early QP subproblems, these phases often require comparable computational effort.
In contrast, dual QP algorithms, such as the Goldfarb-Idnani algorithm 
\cite{goldfarb1983numerically}, reduces computational time by solving 
a dual problem that does not require an explicit feasibility phase.
Several other QP solution strategies have also proven successful
in various SQP implementations.
For instance, SLSQP reformulates QP subproblems as equivalent 
linearly constrained least-squares problems, following the method 
proposed by Schittkowski \cite{schittkowski1982nonlinear2}.
The resulting least-squares problems can be solved efficiently
using methods described by Lawson and Hanson \cite{lawson1995solving}.
Similarly, many SQP methods solve QP subproblems using IP solvers 
\cite{boggs1999global,boggs1999practical,sargent2001new}.
Nonetheless, primal active-set methods are more advantageous as the
problem size increases, owing to their 
superior warm-starting capabilities with accurate active-set estimates.

Using exact Hessians can introduce nonconvexity in QP subproblems,
rendering them NP-hard.
To maintain convexity, a common strategy is to employ positive-definite 
quasi-Newton Hessian approximations.
However, several methods have been developed that can effectively
leverage exact Hessians.
Convexification techniques \cite{gill2015performance,betts1994sparse}
construct a local convex approximation of the QP subproblem using the exact Hessian.
Hybrid approaches \cite{morales2012sequential} first solve the QP subproblem
with a quasi-Newton Hessian, then refine the solution by solving an
equality-constrained QP with the exact Hessian on the identified active set. 
Alternatively, interior point methods can be used to directly solve QPs 
involving exact Hessians, taking advantage of their ability 
to handle indefinite matrices through regularization.

SQP methods designed for large-scale problems employ sparse matrix
techniques to ensure scalability.
Reduced-Hessian approaches, such as those presented in 
\cite{gill2005snopt,tjoa1991simultaneous,eldersveld1992large},
are particularly effective for large-scale problems with moderate 
numbers of degrees of freedom (e.g., up to 2000).
When addressing very large problems with significantly more degrees of freedom, 
QP solvers based on sparse factorizations of the full KKT matrix 
are generally more efficient \cite{betts1994sparse}, especially 
when the problem structure allows for efficient sparse linear algebra operations.

Most globalization strategies rely on
line-search or trust-region methods, in which
candidate steps are accepted based on sufficient decrease
in a merit function or dominance with respect to a filter.
Among merit functions, the $l_1$-penalty function is the most widely used
due to its simplicity; however, it is susceptible to the Maratos effect.
Sophisticated merit functions, such as 
the augmented Lagrangian merit function \cite{gill1986some},
are more robust and do not suffer from this issue.
FilterSQP \cite{fletcher1998user,fletcher2002nonlinear} is a successful 
large-scale, reduced-Hessian method that employs a filter-based strategy
to promote global convergence.

Several SQP algorithms target large-scale equality-constrained problems.
A notable example is the algorithm by Lalee, Nocedal, 
and Platenga \cite{lalee1998implementation}, which builds on the 
trust-region method proposed by Byrd \cite{byrd1987robust} 
and Omojokun \cite{omojokun1989trust}.
This algorithm can utilize second derivative information when available,
but also supports limited-memory BFGS approximations \cite{liu1989limited}.
It forms the basis for the \textit{trust-constr} solver in SciPy.
Another notable equality-constrained method is the 
reduced-Hessian algorithm by Biegler, Nocedal, and Schmid \cite{biegler1995reduced}, 
capable of achieving one-step Q-superlinear convergence with quasi-Newton
approximations of the reduced Hessian.

\section{Methodology}
\label{sec:meth}
In this section, we present the key components of the default SQP algorithm 
implemented within our modular framework.
We also provide an outline of the algorithm and 
discuss practical enhancements to improve robustness and efficiency.

\subsection{The Hessian approximation}
\label{sec:bfgs}
SQP algorithms typically employ positive-definite quasi-Newton 
approximations of the Lagrangian Hessian to ensure 
that the QP subproblems are convex and efficiently solvable,
and to guarantee a descent direction for the merit function during the line search.
In many applications, computing exact second derivatives is difficult or 
computationally intractable, which also necessitates reliance on approximations.
The quasi-Newton BFGS method maintains
good convergence behavior for problems with several
thousand degrees of freedom \cite{gill2015performance}.

In our algorithm, we use the BFGS update,
\begin{equation}
    \widehat{H}_{k+1} = \widehat{H}_k 
    - \frac{1}{d_k^T\widehat{H}_kd_k} \widehat{H}_k d_k d_k^T \widehat{H}_k
    + \frac{1}{w_k^T d_k} w_k w_k^T,
    \label{eq:bfgs}
\end{equation}
where $\widehat{H}_k$ is the current Hessian approximation, 
$d_k = x_{k+1} - x_k$ is the most recent primal step, and 
$w_k = \nabla_x \Lag(x_{k+1}, \lambda_{k+1}) - \nabla_x \Lag(x_k, \lambda_{k+1})$
is the corresponding difference in the Lagrangian gradients
evaluated using the latest multipliers $\lambda_{k+1}$.
The initial approximation $\widehat{H}_0$ is set as the identity matrix.

The BFGS update preserves the positive definiteness of $\widehat{H}_k$ if $w_k^T d_k > 0$.
However, since the Lagrangian Hessian is not necessarily positive-definite 
(even near a local solution),
some iterations may result in $w_k^T d_k \leq 0$.
To address this issue, we adopt the damping scheme proposed by 
Powell \cite{powell1978algorithms},
redefining $w_k$ in (\ref{eq:bfgs}) as
\begin{align}
    \widehat{w_k} = \theta_k w_k + (1-\theta_k)\widehat{H_k} d_k,
\end{align}
where $\theta_k \in [0,1]$ is chosen as
\begin{equation}
    \theta_k = 
    \begin{cases}
        1 & \text{if } w_k^T d_k \geq 0.2\, d_k^T \widehat{H_k} d_k, \\
        \frac{0.8\, d_k^T \widehat{H_k} d_k}{d_k^T \widehat{H_k} d_k - d_k^T w_k} & \text{otherwise}.
    \end{cases}
\end{equation}
This scheme prevents the smallest eigenvalue from decaying 
by more than a factor of $0.2$, 
thereby ensuring hereditary positive-definiteness of $\widehat{H_k}$ and 
descent directions for the merit function, even
when the steps pass through low-curvature regions.

\subsection{The QP solver}
\label{sec:qp}
We solve each QP subproblem using the Python solver \textit{quadprog},
which implements the dual active-set algorithm 
proposed by Goldfarb and Idnani \cite{goldfarb1983numerically} for strictly convex
quadratic programs.
As a dense QP solver, \textit{quadprog} is effective for moderate-scale problems.
However, to maintain computational efficiency when scaling to large-scale problems,
sparse QP solvers would be necessary.

The feasibility phase in a two-phase primal active-set algorithm
can be as expensive as the optimality phase when a feasible point 
is not available a priori---a notable disadvantage during early major iterations 
where the active set changes frequently.
Dual methods offer a computational advantage in this context:
they can be initialized with a dual feasible point, typically obtained 
from the solution of the corresponding unconstrained QP subproblem.
A limitation of this approach is its inability to incorporate warm-starting 
with prior active-set estimates, since primal feasibility from a previous 
QP solution does not necessarily imply dual feasibility for the current subproblem.
For this reason, a switch to a primal active-set algorithm in later iterations 
may be beneficial when reliable active-set estimates become available.

\subsection{Handling infeasible QP}
\label{sec:incons}
Even if the original problem (\ref{eq:np}) is feasible,
constraint linearizations may yield infeasible QP subproblems.
Several strategies exist to address this inconsistency. 
SNOPT, for example, enters a nonlinear elastic mode where it minimizes 
a modified objective that penalizes the $l_1$-norm of the constraint violations.
This modification is equivalent to Tone's second modified QP subproblem 
\cite{tone1982revisions}
and shares strong conceptual ties with 
Fletcher's Sl$_1$QP algorithm \cite{fletcher1985sl1qp}
and Conn's penalty method \cite{conn1973constrained}.
For a broader survey of infeasible linearization handling techniques, 
see Spellucci \cite{spellucci1998new}.

We adopt a modified variant of Powell's simpler strategy \cite{powell1977fast}, 
also implemented in SLSQP.
If a QP subproblem (\ref{eq:qp}) proves to be infeasible,
we instead solve an augmented QP subproblem by introducing an auxiliary variable $\eta$.
The augmented subproblem is given by
\begin{equation}
\begin{split}
    \minimize_{p_k\in\Real^n, \eta\in\Real} & \quad f(x_k) + g(x_k)^T p_k + 
            \frac{1}{2} p_k^T H(x_k,\lambda_k) p_k + \frac{1}{2} \gamma_k \eta^2\\
    \text{subject to} & \quad c_i(x_k)(1-\sigma_i\eta) + \nabla c_i(x_k)^Tp_k \geq 0, \quad i=1,\dots,m, \\
    & \quad 0 \leq \eta \leq 1,
\end{split}
\label{eq:aqp}
\tag{AQP}
\end{equation}
where $\sigma_i$ is defined as
\begin{equation}
    \sigma_i = 
    \begin{cases}
        1 & \text{if } c_i(x_k) < 0, \\
        0 & \text{otherwise},
    \end{cases}
\end{equation}
and $\gamma_k$ is a penalty parameter that 
drives $\eta$ toward zero,
so that the original QP is only perturbed minimally.

This modified variant was proposed by Schittkowski \cite{schittkowski1982nonlinear2},
which retains the distinct treatment for equality and inequality constraints
as in the original approach.
Our implementation extends this further by converting equality constraints into 
pairs of inequality constraints whenever QP subproblems become inconsistent. 
This transformation expands the feasible region and
provides greater freedom to $p_k$ for objective minimization.
Consequently, our approach demonstrates superior performance on benchmark problems, 
particularly in overconstrained cases where equality 
constraints outnumber decision variables.

In its default configuration, OpenSQP sets $\gamma_k=10^6$
whenever infeasibility is detected following a previously feasible iteration. 
If infeasibility persists for multiple consecutive iterations (25, by default)
with a fixed value of $\gamma_k$, 
we increase $\gamma_k$ by a factor of 10 
until reaching a predefined maximum ($10^{12}$, by default).
To solve the augmented QP subproblem (\ref{eq:aqp}), we employ 
the primal active-set solver HiGHS, accessed via the \textit{qpsolvers} library
\cite{qpsolvers2024}. 
Empirical testing across several QP solvers identified HiGHS 
as the most reliable and computationally efficient choice for this task.

\subsection{Global convergence}
\label{sec:globconv}
To guarantee global convergence from arbitrary starting points,
we implement a line search procedure that minimizes a smooth augmented Lagrangian merit function
\cite{gill1986some} along the direction of the QP solution $(p_k, \widehat \lambda_k)$.
The augmented Lagrangian is defined as
\begin{align}
    \Lag_A(x, \lambda, s; \rho) &= f(x) - \lambda^T (c(x)-s) + 
    \frac{1}{2}\sum_{i=1}^m\rho_i(c_i(x)-s_i)^2,
    \label{eq:al}
\end{align}
where $s\in\Real^m$ denotes the slack variables associated with the constraints,
and $\rho\in\Real^m$ is a vector of penalty parameters
that determines the relative weighting of constraint violations.

At each major iteration $k$, following the solution of the quadratic subproblem, we compute the
step size $\alpha_k$ by solving the one-dimensional, unconstrained minimization problem 
\begin{equation}
    \minimize_{\alpha\in\Real} \quad \phi_{\rho}(\alpha) = \Lag_A(x_k + \alpha p_k, \lambda_k + \alpha q_k, s_k + \alpha r_k; \rho),
\label{eq:lsmf}
\end{equation}
where $q_k = \widehat \lambda_k - \lambda_k$ is the step in the dual variables and 
$r_k = \widehat s_k - s_k$ is the step in the slack variables.
Here, $\widehat s_k = c_k + J_k p_k$ represents the slack variables
corresponding to the linearized constraints in the QP subproblem.

We initialize $\lambda_k, s_k,$ and $\rho$ as zero vectors for the first major iteration.
During subsequent iterations, however, some components of $\rho$
may need to be increased to ensure that the search direction 
yields a sufficient decrease in the merit function.
This adjustment can be performed by solving the linearly constrained least-squares problem
\begin{equation}
\begin{split}
    \minimize_{\rho\in\Real^m} & \quad \|\rho\|_2^2 \\
    \text{subject to} & \quad \phi'_{\rho}(0) = -\frac{1}{2}p_k^T \widehat H_k p_k,
\end{split}
\label{eq:lsp}
\end{equation}
where $\phi'_{\rho}(0)$ denotes the derivative of the merit function 
with respect to $\alpha$ at the initial point $\alpha=0$.
The analytical solution $\rho^*$ to this problem is used to
update the penalty parameter vector for the merit function before each line search.
To promote convergence, we adopt the damping scheme 
\begin{equation}
\begin{split}
    \bar \rho_i = \max({\rho_i^*, \hat \rho_i}), \; \text{where } \;
    \hat \rho_i = 
    \begin{cases}
        \rho_i & \text{if } \rho_i < 4(\rho_i^* + \Delta_{\rho}), \\
        (\rho_i(\rho_i^* + \Delta_{\rho}))^{1/2} & \text{otherwise},
    \end{cases}
\end{split}
\label{eq:rho-damped}
\end{equation}
originally proposed by Gill et al. \cite{gill2005snopt}, allowing the penalty parameters to
decrease in a regulated manner during the initial major iterations.
For the first major iteration ($k=0$), $\Delta_{\rho} = 1$.
Subsequently, it is doubled whenever $\| \rho \|_2$ increases following a sequence 
of iterations in which the penalty norm had consistently decreased, 
or conversely, when $\| \rho \|_2$ decreases after a sequence 
of iterations with consistent increases in the penalty norm.
This scheme reduces the fluctuations in $\rho$ during the major iterations
and restricts oscillations in the components of $\rho$ 
to a finite number of occurrences.

After updating $\rho$, we perform a line search using the 
MINPACK \cite{more1980user} implementation,
which enforces the strong Wolfe conditions on the merit function 
to determine the step length $\alpha_k$.
The primal, dual, and slack variables are then updated as 
\begin{align}
    x_{k+1}  = x_k + \alpha_k p_k,  \quad
    \lambda_{k+1} = \lambda_k + \alpha_k q_k,  \quad
    \quad s_{k+1}  = s_k + \alpha_k r_k.
    \label{eq:osqp-ls}
\end{align}
Before each line search, we reset the slack variables $s_k$ by minimizing 
the merit function with respect to $s$, yielding the componentwise update
\begin{equation}
    s_{k_i} = 
    \begin{cases}
        \text{max}(0, c_{k_i}) & \text{if } \rho_{k_i} = 0, \\
        \text{max}(0, c_{k_i} - \lambda_{k_i}/\rho_{k_i} ) & \text{otherwise}.
    \end{cases}
    \label{eq:slackreset}
\end{equation}
This slack reset is performed prior to solving the least-squares problem 
\eqref{eq:lsp} in each major iteration.

\subsection{Convergence criteria}
\label{sec:convergence}
The algorithm terminates successfully at $(x_k, \lambda_k)$ 
if the KKT conditions (\ref{eq:kkt}) are satisfied within specified tolerances.
Let \texttt{opt\_tol} and \texttt{feas\_tol} denote the user-defined 
optimality and feasibility tolerances, respectively.
To reduce the sensitivity of the convergence criteria to the magnitudes 
of the primal and dual variables, we scale these tolerances as
\begin{equation}
    \tau_f=\texttt{feas\_tol} \cdot (1+\|x_k\|_\infty), \quad
    \tau_o=\texttt{opt\_tol} \cdot (1+\|\lambda_k\|_\infty),
\end{equation}
consistent with the scheme used in SNOPT.
\newline
\textbf{Primal feasibility:} A point $x_k$ is considered feasible 
if all constraints $c_i(x_k)$ satisfy
\begin{equation}
    c_i(x_k) \geq -\tau_f.
\end{equation}
\textbf{Optimality conditions:} The pair $(x_k, \lambda_k)$ is deemed optimal 
if $x_k$ is feasible and
\begin{equation}
    [\lambda_k]_i \geq -\tau_o, \quad
    |c_i(x_k)[\lambda_k]_i| \leq \tau_o, \quad
    \|g(x_k)-J(x_k)^T\lambda_k\|_{\infty} \leq \tau_o.
\end{equation}

\subsection{Practical aspects}
\label{sec:practical}
Despite the damped BFGS update \eqref{eq:bfgs}, 
numerical roundoff errors can cause $\widehat{H_k}$ to become indefinite,
particularly when the algorithm remains in regions of negative curvature for 
several consecutive iterations.
In such cases, we discard the current Hessian approximation 
and reset it to the identity matrix, 
as the immediate previous steps, having been taken along directions 
of low curvature, cannot be reliably reused.

A user-provided starting point $x_0$ may not be feasible for some problems. 
SNOPT mitigates this issue by solving a proximal point problem 
to identify the nearest point to $x_0$ that satisfies the linear constraints. 
Since our implementation distinguishes only between variable bounds 
and general constraints, we initialize our algorithm 
from the projection of $x_0$ onto the feasible region defined by the bounds.
If an SQP algorithm starts from such a bound-feasible point, 
all subsequent iterates remain within the bounds, 
allowing the bound constraints to be excluded from the line search.
However, if model functions are undefined at this projected point---an issue that 
occasionally arises---our algorithm reverts to the original 
user-specified starting point $x_0$. 
This contrasts with SNOPT, which terminates immediately 
if function evaluations fail at the computed proximal point.

During the line search, function or derivative evaluations may fail 
when a full unit step is attempted along the search direction $p_k$.
When this occurs, the line search procedure employs a backtracking mechanism, 
reducing the step size along $p_k$ to recover a valid evaluation point.
Separately, the line search may also fail to identify a point that satisfies 
the strong Wolfe conditions, even if evaluations are successful.
In such scenarios, OpenSQP resorts to an inexact line search 
that utilizes an $l_1$-penalty merit function to guide step selection
and promote progress toward a solution.

\subsection{Outline of the algorithm}
\label{sec:alg}

Starting from an initial point $x_0$, 
the OpenSQP algorithm proceeds through the following main steps.

\begin{enumerate}
    \item Compute the bound-feasible proximal point $x_s$ by projecting
          $x_0$ onto the variable bounds.
    \item Evaluate the functions $f(x_s)$, $c(x_s)$ and 
          their gradients $g(x_s)$, $J(x_s)$. 
          If evaluations fail, attempt evaluations at the original $x_0$.
          If evaluations at $x_0$ also fail, terminate.
          Otherwise, set the iteration counter $k \leftarrow 0$ and
          the current iterate $x_k \leftarrow$ (successful $x_s$ or $x_0$).
    \item Initialize the Hessian approximation $\widehat{H}_k \leftarrow I$ 
          (identity matrix), Lagrange multipliers $\lambda_k \leftarrow 0$, 
          slack variables $s_k \leftarrow 0$, and
          penalty parameters $\rho \leftarrow 0$.
    \item Reset the slack variables according to (\ref{eq:slackreset}).
    \item Form the QP subproblem (\ref{eq:qp}) at $x_k$ using $f(x_k)$, 
          $c(x_k)$, $g(x_k)$, $J(x_k)$, and $\widehat{H}_k$.
          Solve the QP to obtain the search direction $(p_k, q_k)$.
          If $\widehat{H}_k$ is insufficiently positive definite 
          for the QP solver, reset $\widehat{H}_k \leftarrow I$ 
          and reattempt QP solution.
          If the QP is infeasible, formulate and solve the augmented QP (\ref{eq:aqp})
          to compute $(p_k, q_k)$. 
          If the augmented QP is also infeasible, terminate after declaring
          the problem as infeasible.
    \item Evaluate the functions $f$, $c$, $g$, and $J$ for the unit step at $x_k+p_k$.
          If any evaluation fails, backtrack along $p_k$ until 
          a valid point $x_k + \beta p_k$ is found where evaluations are successful or 
          the step size $\beta$ becomes too small.
          If evaluations are successful, set the initial trial step length for 
          the line search to $\beta$.
          If the step size $\beta$ falls below a minimum threshold $\beta_{min}$,
          terminate with a status indicating: 
          ``Unable to make progress around undefined region".
    \item Update the penalty parameters $\rho$ as described in Section \ref{sec:globconv}.
    \item Perform a line search to find a step length $\alpha_k$ 
          that satisfies the strong Wolfe conditions for the augmented Lagrangian
          merit function.
          If this fails, attempt a simpler inexact line search based on an 
          $l_1$ penalty merit function.
    \item Update the primal and dual variables: 
          $x_{k+1} \leftarrow x_k + \alpha_k p_k$,
          $\lambda_{k+1} \leftarrow  \lambda_k + \alpha_k q_k$.
    \item Evaluate $f, c, g,$ and $J$ at $x_{k+1}$. 
         (These are typically obtained from the line search; 
         reevaluate only if necessary.)
    \item Update the Hessian approximation $\widehat{H}_{k}$ as described
          in Section \ref{sec:bfgs}, using $d_k=\alpha_k p_k$ and 
          $w_k=g_{k+1}-g_k - (J_{k+1}-J_k)^T \lambda_{k+1}$.
    \item Check if the convergence criteria (see Section \ref{sec:convergence}) are met. 
          If converged, terminate successfully. 
          Else, set $k \leftarrow k+1$ and go to step 4.
\end{enumerate}

\section{Numerical results}
\label{sec:nr}
This section presents a comprehensive evaluation of OpenSQP against 
established nonlinear optimization algorithms. 
Specifically, we benchmark OpenSQP's performance relative to 
SNOPT (a leading commercial SQP optimizer), 
SLSQP (a widely-used, open-source SQP optimizer), 
IPOPT (a widely-used, open-source interior point optimizer), and 
trust-constr (a transparent interior point optimizer in Python).
The comparison utilizes problems from the standardized CUTEst test collection, 
providing an objective assessment across diverse problems.

\subsection{Experimental setup}
The computational experiments were conducted on a MacBook Pro 
with a 2.4 GHz quad-core Intel Core i5 processor and 8 GB of RAM. 
To facilitate a consistent testing environment, all solvers were sourced through 
the modOpt library, which offers a unified interface to all solvers considered
and seamless integration with the CUTEst problem collection via PyCUTEst.
For SLSQP specifically, we utilized the PySLSQP implementation \cite{joshy2024pyslsqp},
which provides more comprehensive optimizer data than SciPy's native implementation.
SNOPT was accessed via the SNOPT-C interface, which 
invokes function and gradient evaluations together.
The specific software versions utilized in this study are as follows: 
\texttt{SNOPT 7.7.7}, \texttt{PySLSQP 0.1.3}, 
\texttt{IPOPT 3.14.11} (with \texttt{MUMPS 5.4.1} linear solver), 
\texttt{SciPy 1.13.1} (for trust-constr),
\texttt{PyCUTEst 1.7.1}, and
QP solvers \texttt{HiGHS 1.8.1} and \texttt{quadprog 0.1.12} for solving
QP subproblems in OpenSQP.

\subsection{Problem selection from CUTEst}
\texttt{PyCUTEst 1.7.1} collection contains $1315$ test problems
with diverse characteristics---ranging from single-variable unconstrained problems 
such as \texttt{BQP1VAR} to large-scale instances such as 
\texttt{BA-L52} with $694,346$ variables and $192,627$ constraints.
While PyCUTEst allows the number of variables and constraints to be
modified for a majority of problems, our tests retain the default dimensions 
specified for each problem instance.
As the goal of our study is to investigate the robustness and reliability of OpenSQP
rather than its computational scalability,
we restrict our study to problems with no more than $100$ variables 
and $100$ constraints.
This restriction also reduces potential issues like memory leaks or 
segmentation faults inherent in extensive test runs with PyCUTEst.
Out of the total $1,315$ problems, $582$ matched our selection criteria. 
Following the exclusion of seven problems due to import incompatibilities 
(namely, 
\texttt{DMN15102LS, DMN15103LS, DMN15332LS, DMN15333LS, DMN37142LS, DMN37143LS,} 
and \texttt{BLEACHNG}), the final test set consisted of $575$ problems.

\subsection{Solver configuration}

To ensure equitable comparisons, we update several default solver parameters. 
A uniform maximum iteration limit of $250$ was imposed on all solvers, 
deemed sufficient for the selected problem dimensions ($m, n \leq 100$). 
IPOPT and SNOPT apply limited-memory Hessian approximations by default 
under certain conditions. 
Specifically, IPOPT always uses a limited-memory update if no Hessian is provided, 
and SNOPT activates limited-memory mode for problems with more than $75$ variables.
To maintain consistency with the full-memory BFGS updates used in the other solvers, 
we increased the limited memory history size to $1000$ for both IPOPT and 
SNOPT, effectively mimicking a full-memory behavior within the 250-iteration limit.
Convergence tolerances were also carefully calibrated. 
For SNOPT and IPOPT, we adopted tolerances from the benchmarking study by 
Gill et al. \cite{gill2015performance}. 
Since OpenSQP employs similar theoretical convergence criteria as SNOPT, 
its tolerances were set identical to those used by SNOPT.
For SLSQP and trust-constr, lacking comprehensive benchmarking literature,
tolerances were tuned empirically across several runs to produce solutions comparable 
in quality to the other solvers. 
Table~\ref{tab:options} provides the complete list of modified solver options.

\begin{table}[h]
\centering
\caption{Non-default solver options}
\begin{tabular}{lll}
\toprule
\textbf{Algorithm} & \textbf{Option} & \textbf{Value} \\
\midrule
\multirow{4}{*}{IPOPT} & \texttt{max\_iter} & \texttt{250} \\
 & \texttt{tol} & \texttt{1.00E-06} \\
 & \texttt{hessian\_approximation} & \texttt{limited-memory} \\
 & \texttt{limited\_memory\_max\_history} & \texttt{1000} \\
\addlinespace
\multirow{3}{*}{OpenSQP} & \texttt{maxiter} & \texttt{250} \\
 & \texttt{opt\_tol} & \texttt{1.22E-04} \\
 & \texttt{feas\_tol} & \texttt{2.00E-06} \\
\addlinespace
\multirow{2}{*}{PySLSQP} & \texttt{maxiter} & \texttt{250} \\
 & \texttt{acc} & \texttt{1.00E-06} \\
\addlinespace
\multirow{4}{*}{SNOPT} & \texttt{Major iterations} & \texttt{250} \\
 & \texttt{Major optimality} & \texttt{1.22E-04} \\
 & \texttt{Major feasibility} & \texttt{2.00E-06} \\
 & \texttt{Hessian updates} & \texttt{1000} \\
\addlinespace
\multirow{3}{*}{trust-constr} & \texttt{maxiter} & \texttt{250} \\
 & \texttt{gtol} & \texttt{2.00E-05} \\
 & \texttt{xtol} & \texttt{2.00E-100} \\
\bottomrule
\label{tab:options}
\end{tabular}
\end{table}

\subsection{Success rates}

We evaluate solver success using algorithm-specific exit conditions:
\begin{enumerate}
    \item SNOPT: exit with the informational message 
    ``\texttt{optimality conditions satisfied}".
    \item IPOPT: exit with the informational message ``\texttt{Optimal Solution Found}".
    \item OpenSQP, PySLSQP, and trust-constr: return \texttt{results[`success']=True}.
\end{enumerate}
Unbounded or infeasible problems are counted as failures for all solvers.

Table \ref{tab:summary} encapsulates the overall success rates, 
providing a measure of algorithm robustness and reliability. 
OpenSQP and SNOPT achieve the highest success rate of $83.30\%$, 
resolving $479$ out of $575$ problems. 
IPOPT follows closely with 475 successful solves ($82.6\%$). 
This study also provides valuable insights as a comprehensive benchmark 
including SLSQP and trust-constr. 
Despite its simpler approach, SLSQP's strong performance, solving over $80\%$ 
of the problems, underscores its continued relevance as a widely-used optimizer. 
Meanwhile, trust-constr lagged behind with a markedly lower 
success rate of $65.7\%$, highlighting limitations in its current implementation 
and indicating room for further improvement.

\begin{table}[htbp]
\caption{Benchmarking summary}
\begin{center}      
    \begin{tabular}{lc}
    \toprule
        \textbf{Algorithm} & \textbf{Successful problems} \\
    \midrule
        IPOPT & 475 (82.61\%) \\
        OpenSQP & 479 (83.30\%) \\
        PySLSQP & 463 (80.52\%) \\
        SNOPT & 479 (83.30\%) \\
        trust-constr & 378 (65.74\%) \\
    \bottomrule
    \end{tabular}
\label{tab:summary}
\end{center}
\end{table}

These findings highlight the robustness and versatility of OpenSQP. 
Crucially, when compared to trust-constr, the only other transparent 
algorithm in this comparison, OpenSQP solved significantly more problems.
This places OpenSQP in a unique position as a transparent solver 
that also delivers performance comparable to leading proprietary 
and well-established open-source solvers.
A key differentiator in success rates among the algorithms is their 
handling of problems with more equality constraints than decision variables.
Only SNOPT and OpenSQP, through their mechanisms for addressing infeasible subproblems, 
could manage such cases.
While IPOPT, SLSQP, and trust-constr failed on all $43$ instances in this category
of overdetermined problems, 
SNOPT and OpenSQP succeeded in solving approximately $20\%$ of these problems.

\subsection{Computational efficiency}

Having established the robustness of each solver, we now assess
the relative efficiency of the solvers using performance and data profiles.
Following Dolan and Mor\'e \cite{dolan2002benchmarking}, 
we construct performance profiles by comparing each solver’s time 
to solve a problem relative to the best-performing solver.
For a set $\mathcal{S}$ of $n_s$ solvers applied to a test set 
$\mathcal{P}$ of $n_p$ problems, let $t_{ps}$ be the performance measure
(e.g., solve time or number of function evaluations) of solver $s\in\mathcal{S}$ 
on problem $p\in\mathcal{P}$.
The performance ratio of solver $s$ on problem $p$ is computed as
\begin{align}
    r_{ps} = 
    \begin{cases}
        \frac{t_{ps}}{\min \{t_{pq}: q\in\mathcal{S}, q \text{ successful on } p\}} & \text{if $s$ was successful}, \\
        \beta r_{max} & \text{otherwise},
    \end{cases}
\end{align}
where $r_{max} = \text{max}\{r_{ps}: p\in\mathcal{P}, s\in\mathcal{S}, 
s \text{ successful on } p\}$ is the maximum performance ratio among all successful
solves over all problems, and $\beta$ is some scalar greater than $1$.
The log-scaled performance function for solver $s$ is then defined as
\begin{align}
    P_s(\tau) = 
        \frac{1}{n_p}| \{p\in\mathcal{P}: \log_2 (r_{ps})\leq \tau\} |,
\end{align}
where $\tau \in [0, \log_2 (r_{max})]$ and $|\cdot|$ denotes the cardinality of a set.
The value of $P_s(\tau)$ represents the proportion of problems that was solved by
solver $s$ within $2^{\tau}$ times the best solve time or number of function evaluations.
Performance profiles refer to the plots of performance functions
with the performance measure as solve time.
Data profiles, as introduced in \cite{more2009benchmarking}, 
use the total number of function evaluations as the performance measure 
instead of solve time.

\begin{figure}[ht]
    \centering
    \includegraphics[width=0.45\linewidth]{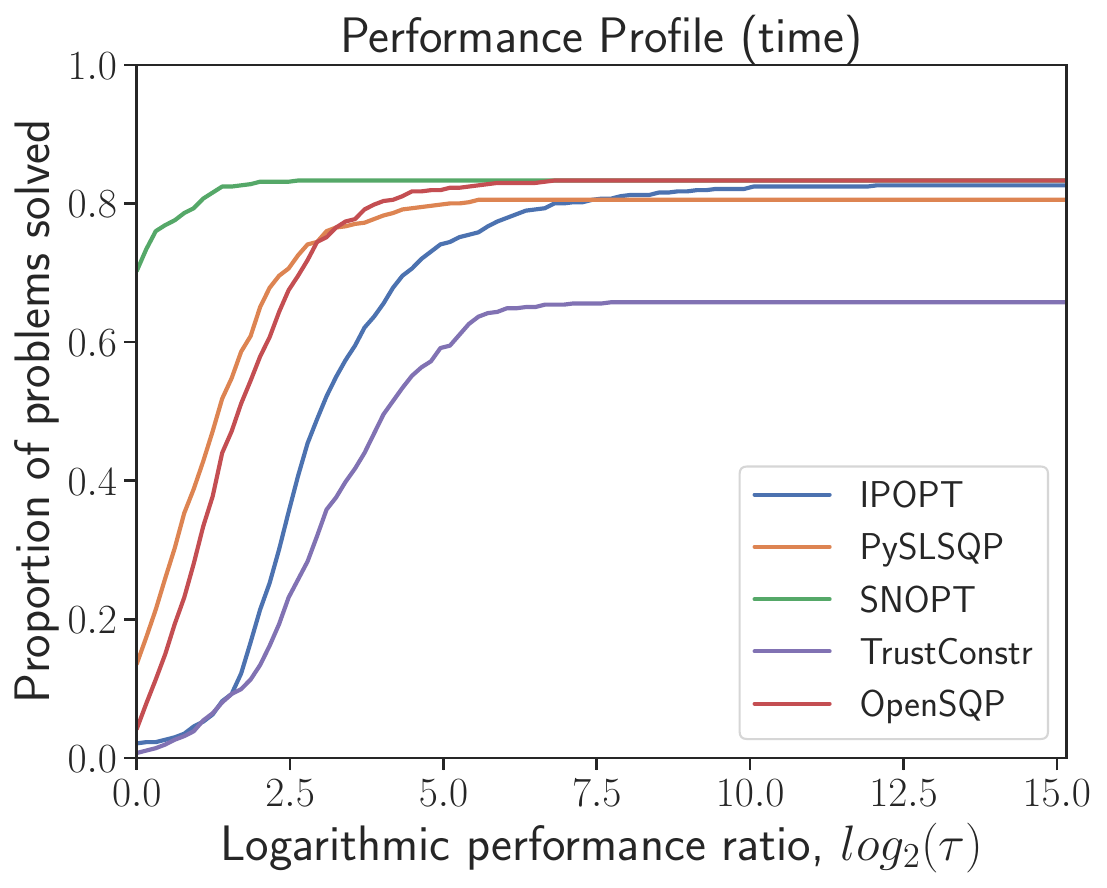}
    \includegraphics[width=0.45\linewidth]{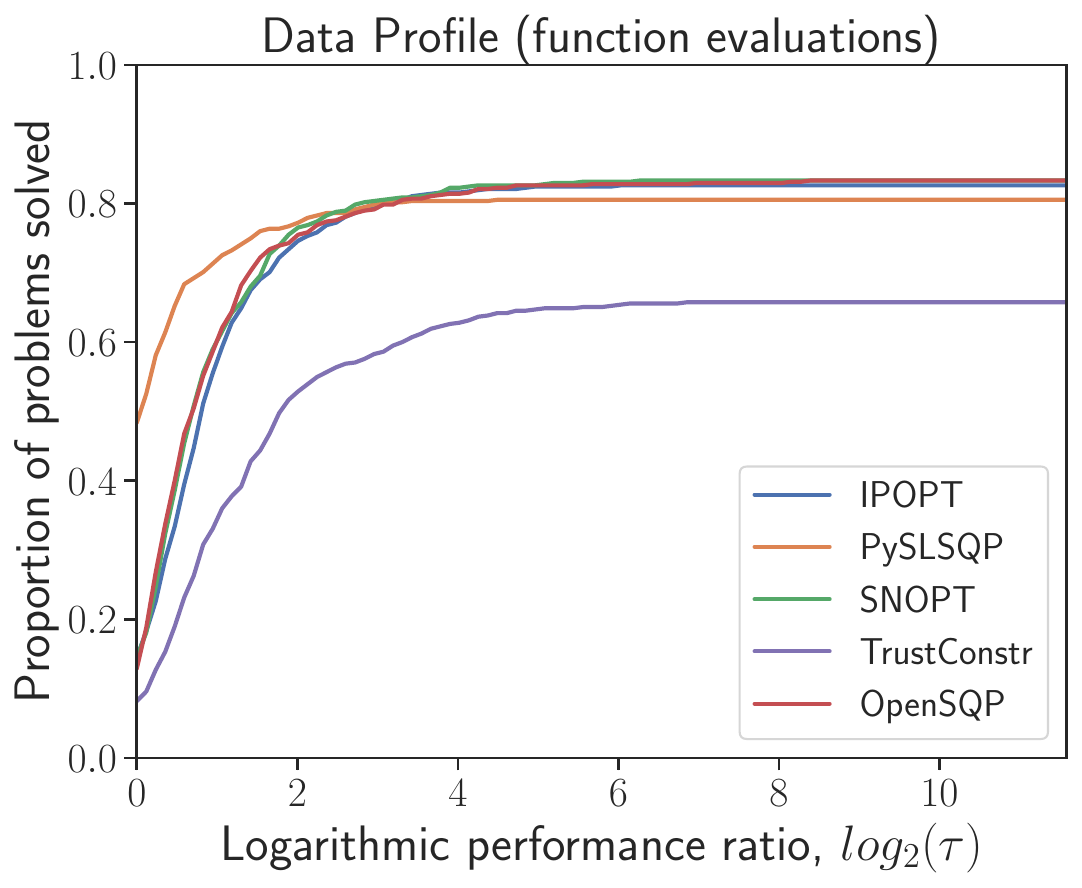}
    \caption{Performance and data profiles of IPOPT, PySLSQP, OpenSQP, SNOPT, and trust-constr on 575 problems from the CUTEst collection.}
    \label{fig:cutest}
\end{figure}

Figure \ref{fig:cutest} presents the performance and data profiles of all solvers 
for the $575$ selected problems.
For data profiles, we define ``function evaluations'' as the combined number 
of objective, objective gradient, constraint, and 
constraint Jacobian evaluations.
Several key observations emerge from these profiles.
\begin{enumerate}
\item \textbf{Competitiveness of OpenSQP}: In terms of function evaluations, 
OpenSQP performs on par with the state-of-the-art algorithms IPOPT and SNOPT.
Despite being implemented in Python,
OpenSQP's performance in time approaches that of the Fortran-based SLSQP.
\item \textbf{Speed of SNOPT}: SNOPT consistently achieved the fastest solve times
of all solvers. 
This likely stems from its efficient and highly optimized linear algebra routines.
\item \textbf{Slower performance of IPOPT}: 
IPOPT appears less time-efficient in this setting. 
It is important to acknowledge, however, that for problems of modest scale, 
computational times may reflect interface overheads more than algorithmic efficiency, 
making data profiles particularly valuable for a fair comparison.
\item \textbf{Efficiency of SLSQP}: Despite solving fewer problems than the top
three solvers, 
SLSQP demonstrates the highest efficiency in terms of function evaluations 
for the problems it solves successfully.
\item \textbf{Limitations of trust-constr}:
Both performance and data profiles consistently indicate that trust-constr is less 
robust and generally less efficient than the other four contenders in this benchmark,
suggesting that further tuning or implementation improvements 
may be necessary for broader applicability.

\end{enumerate}

Overall, these results demonstrate that OpenSQP offers a strong combination of 
transparency, flexibility, and performance, achieving results comparable to mature, 
production-level solvers while retaining the advantages of a modular, 
reconfigurable design.
This positions OpenSQP as a compelling alternative to established optimization
software.

\section{Conclusions}
\label{sec:conclusion}

We presented OpenSQP, a modular and reconfigurable SQP algorithm designed to
support optimization practitioners and researchers
in developing SQP algorithms tailored to specific needs.
Implemented in Python, its software design 
enables users to easily modify key components, such as merit functions, 
line search algorithms, Hessian approximations, and QP solvers, to better 
suit their applications.
Because of its explicit support for modular and transparent optimizer development, 
we develop and distribute OpenSQP within the modOpt optimization library.
The algorithm and its documentation are hosted open-source 
on modOpt's GitHub repository.

In Section \ref{sec:bg}, we presented some background on nonlinear programming,
followed by a brief review of sequential quadratic programming.
Section \ref{sec:meth} discussed the algorithmic details of 
the default configuration of OpenSQP, which employs
a smooth augmented Lagrangian merit function for the line search and
a positive-definite BFGS Hessian approximation.
Section \ref{sec:nr} assessed the performance of the
default configuration of OpenSQP using a large set of problems 
from the CUTEst test set.
The benchmark results demonstrated reliability and efficiency comparable to 
state-of-the-art NLP algorithms such as SNOPT and IPOPT.

Sequential quadratic programming represents a broad area of research.
The algorithm presented in this paper serves as an accessible entry point
for SQP researchers and practitioners.
Several improvements could further enhance its capabilities,  
including incorporating a limited-memory Hessian approximation 
for better scalability,
modifying the algorithm to utilize exact Hessians or Hessian-vector 
products for improved performance, 
and hybridizing the QP solver by combining interior-point and active-set methods 
for efficiently solving QP subproblems in large-scale problems.


%

\section*{Funding}
This material is based upon work supported by 
the National Science Foundation under Grant No. 1917142 and 
by NASA under Award No. 80NSSC23M0217.

\section*{Declarations}

\section*{Conflict of interest}
On behalf of all authors, the corresponding author states 
that there is no conflict of interest.

%

\section*{Code availability}

The OpenSQP solver is available open-source as part of the modOpt
Python library.
The source-code can be accessed through modOpt's GitHub repository at
\href{https://github.com/lsdolab/modopt/}{https://github.com/lsdolab/modopt/}.
The code for all numerical studies and benchmarks presented 
in this paper is available in the file 
\href{https://github.com/lsdolab/modopt/tree/main/examples/ex_19_opensqp_benchmark.py}{\texttt{examples/ex\_19\_opensqp\_benchmark.py}} within the same repository.
Comprehensive documentation, including installation instructions, usage guidelines, 
and API references, can be found at 
\href{https://modopt.readthedocs.io/}{https://modopt.readthedocs.io/}.

\bibliographystyle{spmpsci}      
\bibliography{references}   

\end{document}